\documentclass{gtart}

\def\ifplaintex{\expandafter\ifx\csname documentclass\endcsname\relax}

\def\gtp{{\mathsurround=0pt\it $\cal G\mskip-2mu$eometry \&\ 
$\cal T\!\!$opology $\cal P\!$ublications}}  

\def\recd{{\small Received:\qua\receiveddate\ifx\reviseddate\relax
\else\qquad Revised:\qua\reviseddate\fi\par}} 

\def\lognumber#1{\def\thelognumber{#1}}
\def\volumenumber#1{\def\thevolumenumber{#1}}
\def\volumeyear#1{\def\thevolumeyear{#1}}
\def\papernumber#1{\def\thepapernumber{#1}}
\def\pagenumbers#1#2{\def\startpage{#1}\def\finishpage{#2}}
\def\published#1{\def\publishdate{#1}}

\def\received#1{\def\receiveddate{#1}}
\def\revised#1{\def\reviseddate{#1}}
\def\accepted#1{\def\accepteddate{#1}}

\let\\\par\let\thelognumber\relax\let\thevolumenumber\relax
\let\thepapernumber\relax\let\thevolumeyear\relax\let\startpage\relax
\let\finishpage\relax\let\publishdate\relax\let\receiveddate\relax
\let\reviseddate\relax\let\accepteddate\relax\let\theasciititle\relax
\let\theasciiauthors\relax
\let\theasciiabstract\relax

\let\theasciiemail\relax

\ifplaintex
\font\logobig=cmssbx10 scaled 3836
\font\logomed=cmssbx10 scaled 2557
\else
\font\logobig=cmssbx10 scaled 4200
\font\logomed=cmssbx10 scaled 2800
\fi

\long\def\makeagttitle{   
\count0=\startpage
\agt\hfill      
\hbox to 45truept{\vbox to 0pt{\vglue -13truept{\logomed A\kern -.37em{\logobig 
T}\kern -.38em G}\vss}\hss}
\break
{\small Volume \thevolumenumber\ (\thevolumeyear)
\startpage--\finishpage\nl
Published: \publishdate}

\vglue .25truein

{\parskip=0pt\leftskip 0pt plus
1fil\def\\{\par\smallskip}{\Large\bf\thetitle}\par\medskip} \vglue
0.05truein

{\parskip=0pt\leftskip 0pt plus 1fil\def\\{\par}{\sc\theauthors}
\par\medskip}%
 
\vglue 0.03truein 

{\small\leftskip 25truept\rightskip 25truept{\bf Abstract}\stdspace\theabstract

{\bf AMS Classification}\stdspace\theprimaryclass
\ifx\thesecondaryclass\relax\else; \thesecondaryclass\fi\par
{\bf Keywords}\stdspace \thekeywords\par}\vglue 7truept

}   

\ifplaintex
\font\phead=cmsl9 scaled 950
\font\pnum=cmbx10 scaled 913
\font\pfoot=cmsl9 scaled 950
\headline{\vbox to 0pt{\vskip -4.5mm\line{\small\phead\ifnum
\count0=\startpage ISSN 1472-2739 (on-line) 1472-2747 (printed)
\hfill {\pnum\folio}\else\ifodd\count0\def\\{ }%
\ifx\theshorttitle\relax\thetitle\else\theshorttitle\fi\hfill{\pnum\folio}
\else\def\\{ and }{\pnum\folio}\hfill\ifx\theshortauthors\relax\theauthors
\else\theshortauthors\fi\fi\fi}\vss}}
\footline{\vbox to 0pt{\vglue 0mm\line{\small\pfoot\ifnum\count0=\startpage
\copyright\ \gtp\hfill\else
\agt, Volume \thevolumenumber\ (\thevolumeyear)\hfill\fi}\vss}}
\else
\font=cmsl9 scaled 1050
\font\lnum=cmbx10 
\font=cmsl9 scaled 1050
\makeatletter
\def\@oddhead{{\small\lhead\ifnum\count0=\startpage ISSN 1472-2739 
(on-line) 1472-2747 (printed)\hfill {\lnum\number\count0}\else\ifodd\count0
\def\\{ }\ifx\theshorttitle\relax \thetitle \else\theshorttitle\fi\hfill
{\lnum\number\count0}\else\def\\{ and }{\lnum\number\count0}
\hfill\ifx\theshortauthors\relax 
\theauthors\else\theshortauthors\fi\fi\fi}}\def\@evenhead{\@oddhead}
\def\@oddfoot{\small\lfoot\ifnum\count0=\startpage\copyright\ \gtp\hfill\else
\agt, Volume \thevolumenumber\ (\thevolumeyear)\hfill\fi}
\def\@evenfoot{\@oddfoot}
\makeatother
\fi
\let\maketitlepage\makeagttitle
\let\makeshorttitle\maketitlepage
\let\maketitle\maketitlepage

\newwrite\gtoutfile
\long\gdef\makeheadfile{  
{\def\\{, }\def\s{ }
\immediate\openout\gtoutfile head.xxx
\immediate\write\gtoutfile{To: math@arxiv.org}
\immediate\write\gtoutfile{Subject: put OR rep NNNNN:ppppp}
\immediate\write\gtoutfile{--text follows this line--}
\immediate\write\gtoutfile{Proxy-for: \ifx\theasciiauthors\relax
\theauthors\else\theasciiauthors\fi\s<\ifx\theasciiemail\relax\theemail\else\theasciiemail\fi>}
\immediate\write\gtoutfile{\noexpand\\}
\immediate\write\gtoutfile{Authors: \ifx\theasciiauthors\relax
\theauthors\else\theasciiauthors\fi}
{\def\\{ }\immediate\write\gtoutfile{Title: \ifx\theasciititle\relax
\thetitle\else\theasciititle\fi}}
\immediate\write\gtoutfile{Subj-class: GT or SG, GR etc}
\immediate\write\gtoutfile{MSC-class: \theprimaryclass\ifx\thesecondaryclass\relax\else, \thesecondaryclass\fi}
\immediate\write\gtoutfile{Journal-ref: Algebr. Geom. Topol. \thevolumenumber\s
(\thevolumeyear) \startpage-\finishpage}
\immediate\write\gtoutfile{Comments: Published by Algebraic and
Geometric Topology at}
\immediate\write\gtoutfile{\s\s\s  http://www.maths.warwick.ac.uk/agt/AGTVol\thevolumenumber/agt-\thevolumenumber-\thepapernumber.abs.html}
\immediate\write\gtoutfile{\noexpand\\}
\immediate\write\gtoutfile{}
\ifx\theasciiabstract\relax
\immediate\write\gtoutfile{\theabstract}\else
\immediate\write\gtoutfile{\theasciiabstract}\fi
\immediate\write\gtoutfile{}
\immediate\write\gtoutfile{\noexpand\\}
\immediate\write\gtoutfile{}
\immediate\closeout\gtoutfile}}  

\def\maketitlepage{\makeagttitle\makeheadfile}
\let\makeshorttitle\maketitlepage
\let\maketitle\maketitlepage

\def\ifplaintex{\expandafter\ifx\csname documentclass\endcsname\relax}

\def\gtp{{\mathsurround=0pt\it $\cal G\mskip-2mu$eometry \&\ 
$\cal T\!\!$opology $\cal P\!$ublications}}  

\def\recd{{\small Received:\qua\receiveddate\ifx\reviseddate\relax
\else\qquad Revised:\qua\reviseddate\fi\par}} 

\def\lognumber#1{\def\thelognumber{#1}}
\def\volumenumber#1{\def\thevolumenumber{#1}}
\def\volumeyear#1{\def\thevolumeyear{#1}}
\def\papernumber#1{\def\thepapernumber{#1}}
\def\pagenumbers#1#2{\def\startpage{#1}\def\finishpage{#2}}
\def\published#1{\def\publishdate{#1}}

\def\received#1{\def\receiveddate{#1}}
\def\revised#1{\def\reviseddate{#1}}
\def\accepted#1{\def\accepteddate{#1}}

\let\\\par\let\thelognumber\relax\let\thevolumenumber\relax
\let\thepapernumber\relax\let\thevolumeyear\relax\let\startpage\relax
\let\finishpage\relax\let\publishdate\relax\let\receiveddate\relax
\let\reviseddate\relax\let\accepteddate\relax\let\theasciititle\relax
\let\theasciiauthors\relax
\let\theasciiabstract\relax

\let\theasciiemail\relax

\ifplaintex
\font\logobig=cmssbx10 scaled 3836
\font\logomed=cmssbx10 scaled 2557
\else
\font\logobig=cmssbx10 scaled 4200
\font\logomed=cmssbx10 scaled 2800
\fi

\long\def\makeagttitle{   
\count0=\startpage
\agt\hfill      
\hbox to 45truept{\vbox to 0pt{\vglue -13truept{\logomed A\kern -.37em{\logobig 
T}\kern -.38em G}\vss}\hss}
\break
{\small Volume \thevolumenumber\ (\thevolumeyear)
\startpage--\finishpage\nl
Published: \publishdate}

\vglue .25truein

{\parskip=0pt\leftskip 0pt plus
1fil\def\\{\par\smallskip}{\Large\bf\thetitle}\par\medskip} \vglue
0.05truein

{\parskip=0pt\leftskip 0pt plus 1fil\def\\{\par}{\sc\theauthors}
\par\medskip}%
 
\vglue 0.03truein 

{\small\leftskip 25truept\rightskip 25truept{\bf Abstract}\stdspace\theabstract

{\bf AMS Classification}\stdspace\theprimaryclass
\ifx\thesecondaryclass\relax\else; \thesecondaryclass\fi\par
{\bf Keywords}\stdspace \thekeywords\par}\vglue 7truept

}   

\ifplaintex
\font\phead=cmsl9 scaled 950
\font\pnum=cmbx10 scaled 913
\font\pfoot=cmsl9 scaled 950
\headline{\vbox to 0pt{\vskip -4.5mm\line{\small\phead\ifnum
\count0=\startpage ISSN 1472-2739 (on-line) 1472-2747 (printed)
\hfill {\pnum\folio}\else\ifodd\count0\def\\{ }%
\ifx\theshorttitle\relax\thetitle\else\theshorttitle\fi\hfill{\pnum\folio}
\else\def\\{ and }{\pnum\folio}\hfill\ifx\theshortauthors\relax\theauthors
\else\theshortauthors\fi\fi\fi}\vss}}
\footline{\vbox to 0pt{\vglue 0mm\line{\small\pfoot\ifnum\count0=\startpage
\copyright\ \gtp\hfill\else
\agt, Volume \thevolumenumber\ (\thevolumeyear)\hfill\fi}\vss}}
\else
\font=cmsl9 scaled 1050
\font\lnum=cmbx10 
\font=cmsl9 scaled 1050
\makeatletter
\def\@oddhead{{\small\lhead\ifnum\count0=\startpage ISSN 1472-2739 
(on-line) 1472-2747 (printed)\hfill {\lnum\number\count0}\else\ifodd\count0
\def\\{ }\ifx\theshorttitle\relax \thetitle \else\theshorttitle\fi\hfill
{\lnum\number\count0}\else\def\\{ and }{\lnum\number\count0}
\hfill\ifx\theshortauthors\relax 
\theauthors\else\theshortauthors\fi\fi\fi}}\def\@evenhead{\@oddhead}
\def\@oddfoot{\small\lfoot\ifnum\count0=\startpage\copyright\ \gtp\hfill\else
\agt, Volume \thevolumenumber\ (\thevolumeyear)\hfill\fi}
\def\@evenfoot{\@oddfoot}
\makeatother
\fi
\let\maketitlepage\makeagttitle
\let\makeshorttitle\maketitlepage
\let\maketitle\maketitlepage

\newwrite\gtoutfile
\long\gdef\makeheadfile{  
{\def\\{, }\def\s{ }
\immediate\openout\gtoutfile head.xxx
\immediate\write\gtoutfile{To: math@arxiv.org}
\immediate\write\gtoutfile{Subject: put OR rep NNNNN:ppppp}
\immediate\write\gtoutfile{--text follows this line--}
\immediate\write\gtoutfile{Proxy-for: \ifx\theasciiauthors\relax
\theauthors\else\theasciiauthors\fi\s<\ifx\theasciiemail\relax\theemail\else\theasciiemail\fi>}
\immediate\write\gtoutfile{\noexpand\\}
\immediate\write\gtoutfile{Authors: \ifx\theasciiauthors\relax
\theauthors\else\theasciiauthors\fi}
{\def\\{ }\immediate\write\gtoutfile{Title: \ifx\theasciititle\relax
\thetitle\else\theasciititle\fi}}
\immediate\write\gtoutfile{Subj-class: GT or SG, GR etc}
\immediate\write\gtoutfile{MSC-class: \theprimaryclass\ifx\thesecondaryclass\relax\else, \thesecondaryclass\fi}
\immediate\write\gtoutfile{Journal-ref: Algebr. Geom. Topol. \thevolumenumber\s
(\thevolumeyear) \startpage-\finishpage}
\immediate\write\gtoutfile{Comments: Published by Algebraic and
Geometric Topology at}
\immediate\write\gtoutfile{\s\s\s  http://www.maths.warwick.ac.uk/agt/AGTVol\thevolumenumber/agt-\thevolumenumber-\thepapernumber.abs.html}
\immediate\write\gtoutfile{\noexpand\\}
\immediate\write\gtoutfile{}
\ifx\theasciiabstract\relax
\immediate\write\gtoutfile{\theabstract}\else
\immediate\write\gtoutfile{\theasciiabstract}\fi
\immediate\write\gtoutfile{}
\immediate\write\gtoutfile{\noexpand\\}
\immediate\write\gtoutfile{}
\immediate\closeout\gtoutfile}}  

\def\maketitlepage{\makeagttitle\makeheadfile}
\let\makeshorttitle\maketitlepage
\let\maketitle\maketitlepage

\lognumber{15}
\volumenumber{1}
\volumeyear{2001}
\papernumber{15}
\pagenumbers{311}{319}
\received{8 February 2001}
\revised{16 May 2001}
\accepted{22 May 2001}
\published{24 May 2001}

\usepackage{amssymb}

\def\f{\varphi}
\def\s{\sigma}
\def\C{{\mathbb C}}

\def\P{{\mathbb P}}

\def\R{{\mathbb R}}
\def\Z{{\mathbb Z}}
\newtheorem{thm}{Theorem}
\newtheorem{prop}[thm]{Proposition}

\newtheorem{lemma}[thm]{Lemma}
\newtheorem{cor}[thm]{Corollary}

\begin{document}
\title{Free group automorphisms, invariant orderings\\and 
topological applications}
\author{Dale Rolfsen\\Bert Wiest}
\address{Mathematics Department, University of British Columbia,
1984 Mathematics Road, Vancouver BC, Canada V6T 1Z2}
\email{rolfsen@math.ubc.ca, bertw@pims.math.ca}
\begin{abstract} We are concerned with orderable groups
and particularly those with orderings invariant not only under
multiplication, but also under a given automorphism or family of
automorphisms.
Several applications
to topology are given: we prove that the fundamental groups of hyperbolic
nonorientable surfaces, and the groups of certain fibred knots are
bi-orderable.
Moreover, we show that the pure braid groups associated with hyperbolic
nonorientable surfaces are left-orderable. \end{abstract}
\primaryclass{6F15}
\secondaryclass{57M05}
\keywords{Ordered group, surface group, knot group, surface braid group}
\maketitle

\section{Introduction}  If $G$ is a group, and $<$ a strict total
ordering of its elements, we say that $(G,<)$ is a {\it left-ordered}
group
if $x<y \Leftrightarrow zx<zy$ for all $x,y,z \in G$, and
{\it bi-ordered} if the ordering is also right-invariant:
$x<y \Leftrightarrow xz<yz$.  Surprisingly many groups are left-orderable
or even bi-orderable; for example free groups are bi-orderable \cite{MR},
although this is by no means obvious (one method is described in Section 3).

If $\f\co G \to G$ is an automorphism, then we say that
an ordering $<$ of $G$ is {\it invariant} under $\f$
(or {\it respected} by $\f$) if
$x<y \Leftrightarrow \f(x) < \f(y)$ for all $x,y \in G$.
As we will see, it may or may not be possible to find such an
invariant ordering, depending on the nature of the $\f$.
In fact all our results on such invariant orderings involve free groups.

The goal of this paper is to establish new orderability results for
several families of groups which arise in topology.  We use a sort of
bootstrap process involving extensions
and orderings of free subgroups which also invariant under certain
automorphisms. We show the following results:

(1)\qua The fundamental groups of all closed surfaces, orientable or not,
are bi-orderable, with the exceptions of the projective plane and the
Klein bottle (cf.\  \cite{BRW}).

(2)\qua The pure braid groups $PB_n(N)$ associated with nonorientable
surfaces $N \ne \P^2$ are left-orderable.

(3)\qua The fundamental groups of certain punctured-torus bundles over the
circle are bi-orderable; examples include the figure-of-eight knot group.

Some of our technical arguments are directly adapted from \cite{G-M}
(which in turn built on ideas from \cite{KR}).
In these papers it was proved that the pure braid groups, and more
generally all pure surface braid groups (in any orientable compact 
surface) are bi-orderable.  Moreover, \cite{G-M} shows that
$PB_n(N)$ is {\it not} bi-orderable for nonorientable surfaces $N$.

More general surveys on orderable groups and their r\^ole in low-dimensional
topology, which place our results in their proper context, can be found in
the introductory sections of \cite{BRW} and \cite{SW}.

We thank Steven Boyer for many helpful discussions; some of the present
results arose directly from our work on the orderability of 3-manifold
groups. Thanks also to Juan Gonzales-Meneses for inspiring discussions.
Dale Rolfsen was partially supported by research grant of
the Canadian Natural Sciences and Engineering Research Council, and
Bert Wiest by a PIMS postdoctoral fellowship.


\section{Orderable groups and extensions}

Groups which are left-orderable are easily seen to be torsion-free.
However, left-orderable groups enjoy certain advantages over merely
torsion-free ones; for instance, they are known to satisfy the (still open)
zero-divisor conjecture, which states that the group ring $\Z G$ of a
torsion free group $G$ should have no zero-divisors.
The group rings of bi-orderable groups are known to have an even
stronger property, due to Malcev and Neumann
(and conjectured to be true for left-orderable groups as well):
if $G$ is a bi-orderable group, then  $\Z G$ embeds in a division algebra.
We refer the reader to \cite{MR} for proofs and more general statements.

It is easy to verify that a group $G$ is left-orderable if and only if
there is a subset $P \subset G$ which does not contain the identity
element, is closed under multiplication and such that for every $g \ne 1$,
exactly one of $g$ and $g^{-1}$ belongs to $P$.  Given such a $P$, define
$g < h \Leftrightarrow g^{-1}h \in P$.  On the other hand, given a
left-ordering $<$, define $P$ to be the positive cone $P = \{g \in G : 1 <
g\}$.  Note that if one instead used the criterion $g < h \Leftrightarrow
h g^{-1} \in P$, a right-invariant ordering would result: a group is
right-orderable if and only if it is left-orderable.

A left-ordering is bi-invariant if and only if its positive cone $P$ is
normal:  $g^{-1}Pg \subset P ~~ \forall g \in G$.  Moreover, the
ordering is invariant under an automorphism $\f$ if and only
$\f(P) \subset P$.
Following is one of the reasons to be interested in orderings invariant
under automorphisms. Its straightforward proof is left to the reader.

\begin{lemma} \label{extensions}\sl
Suppose we have a short exact sequence of groups
$$
1 \longrightarrow F \ \hookrightarrow \ G \stackrel{p}{\longrightarrow}
H:=G/F \longrightarrow 1.
$$
If $F$ and $H$ are left-orderable, then
$G$ is left-orderable, with positive cone $P_G := p^{-1}(P_H) \cup P_F$.
If $F$ and $H$ are bi-ordered, then the same formula defines the
positive cone for a bi-ordering of $G$ if and only if
$g^{-1}P_Fg \subset P_F$ for all $g \in G$, that is, if and only if the
ordering of $F$ is invariant under conjugation by elements of $G$.
\end{lemma}

This gives us a strategy to prove that a group $G$ is bi-orderable.  One
finds a convenient normal subgroup $F$ such that $G/F$ is bi-orderable
and in addition $F$ can be given a bi-ordering invariant under
conjugation by elements of $G$.


\section{Surface groups}

It has been known for a long time \cite{Baumslag, Long} that the
fundamental groups of \emph{orientable} surfaces are bi-orderable
(and in particular left-orderable -- an interesting non-constructive
argument in \cite{bowditch} also shows this).
Our aim is to generalize these results to non-orientable surfaces, using
very different techniques.

We shall denote the connected sum of $n$ projective planes by $n \P^2$.
We recall that $(n+2)\P^2$ is homeomorphic to the connected sum of a torus
and $n$ projective planes. Consider $3\P^2\cong T^2 \# \P^2$, the
nonorientable surface  with Euler characteristic $-1$; this surface
will be the key to our analysis.  Note that $n\P^2$ has a hyperbolic
structure
for $n \geqslant 3$, whereas $\P^2$ is spherical and the Klein bottle $2\P^2$
is Euclidean.

\begin{prop} \label{3P2} \sl
The group $\pi_1(3 \P^2)$ is bi-orderable.  \end{prop}

\begin{thm} \label{surfaces} \sl
If $N$ is any connected surface other than the
projective plane  $\P^2$ or Klein bottle $2\P^2$, then $\pi_1(N)$ is
bi-orderable. For $N=$Klein bottle, $\pi_1(N)$ is left-orderable.
\end{thm}

\begin{proof}[Proof of theorem \ref{surfaces}]
Let us first see how theorem \ref{surfaces} follows from proposition
\ref{3P2}. If $N$ is noncompact, or if $\partial N$
is nonempty, then $\pi_1(N)$ is a free group, and therefore bi-orderable.
Thus we are reduced to considering closed surfaces. According to the
standard classification, such surfaces are either a connected sum of tori,
or of projective planes in the nonorientable case.

We shall first consider non-orientable surfaces.
Of course $\pi_1(\P^2) \cong \Z_2$ is certainly not
left-orderable.   For $N = 2\P^2$, the Klein bottle,
$$\pi_1(N) \cong \langle x, y: xyx^{-1} = y^{-1}\rangle$$
is a well-known example of a group which is left-orderable (being an
extension of $\Z$ by $\Z$), but not bi-orderable, as the defining relation
would lead to a contradiction.

By proposition \ref{3P2}, the surface $M=3\P^2\cong T^2 \# \P^2$ has
bi-orderable fundamental group.
We shall picture it as a torus with a small disk removed, and replaced by
sewing in a M\"obius band. Consider an $n$-fold cover of the torus by
itself, and modify the covering by replacing one disk downstairs, and
$n$ disks upstairs, by M\"obius bands.  This gives a covering of $M$ by
the connected sum of a torus with $n$ copies of $\P^2$.  Thus the
fundamental group of $(n+2)\P^2$ injects in that of $3\P^2$, and is
therefore bi-ordered.

We now turn to orientable closed surfaces. The cases of genus zero or one
being easy, we consider a closed
surface of genus  $g \geqslant 2$.  This surface is the oriented double
cover of  $(g+1)\P^2$.  Therefore its fundamental group is a subgroup of
a bi-orderable group. This completes the proof of theorem \ref{surfaces},
assuming proposition \ref{3P2}. \end{proof}

\begin{cor}\label{PBnLO}\sl The pure braid group $PB_m(N)$ on $m$ strands
in a compact 
surface $N\neq \P^2$ is left-orderable. 
\end{cor}

If $N$ is non-orientable, this is the strongest possible result, because
by \cite{G-M} the groups $PB_m(N)$ are definitely not bi-orderable.
If $N$ is orientable, the result is redundant, since in this case $PB_m(N)$
is known to be even bi-orderable \cite{G-M}. We also remark that it is not
known which non-pure surface braid groups are left-orderable.

\begin{proof}[Proof of corollary \ref{PBnLO}]
We shall proceed by induction on $m$. For $m=1$,
we have that $PB_1(N)=\pi_1(N)$ is left-orderable by theorem
\ref{surfaces}.
Moreover, we have a short exact sequence
$$
1 \longrightarrow \pi_1(N-\{m \hbox{ points}\}) \ \hookrightarrow \
PB_{m+1}(N)
\stackrel{\phi}{\longrightarrow} PB_{m}(N) \longrightarrow 1,
$$
where $\phi$ is induced by forgetting one of the strands. Now
$\pi_1(N-\{m \hbox{ points}\})$ is free and hence left-orderable,
and $PB_{m}(N)$ is left-orderable by induction; an application of
lemma \ref{extensions} completes the  induction step. \end{proof}\\

\begin{proof}[Proof of proposition \ref{3P2}] Let $M=3\P^2 \cong
T^2 \# \P^2$. Our strategy for constructing
a bi-invariant order on $\pi_1(M)$ is to apply
lemma \ref{extensions}, where the normal subgroup $F$ of $G=\pi_1(M)$
will be chosen so that $G/F\cong \Z^2$.

To define the subgroup $F$, we note that $G$ has presentation
$$
G = \langle a, b, c : aba^{-1}b^{-1} = c^2 \rangle,
$$
where $a$ and $b$ represent free generators of the punctured
torus in $M$ and $c$ the generator corresponding to the central curve of
the M\"obius band in $M$. We define
$F = \langle\langle c \rangle\rangle$, the normal subgroup generated by
$c$.  Note that a word in $a,b,c$ belongs to $F$ if and only its exponent
sums in $a$ and $b$ are both zero.

The covering $\widetilde{M}$ of $M$ with $\pi_1(\widetilde{M}) = F$ is
very easy to imagine: consider the universal covering $\R^2 \to T^2$,
and modify $\R^2$ by taking a family of small
disks (say $D_{i,j}$) centered at the integral points $(i,j) \in \Z^2$.
Remove each of these and replace by a M\"obius band $B_{i,j}$.
This defines a covering $\widetilde{M} \to M$.  The group of covering
translations of $\widetilde{M}$ is just $\Z^2$, with $(m,n)$ acting by
translation  $(x,y) \to (x+m,y+n)$ on the $\R^2$ part of $\widetilde{M}$,
and
taking each $B_{i,j}$ to $B_{i+m,j+n}$.  Therefore we have (as required)
an exact sequence
$$
1 \to \pi_1(\widetilde{M}) \to \pi_1(M) \to \Z^2 \to 1.
$$
We now turn to the task of proving that the orderability hypotheses of
lemma \ref{extensions} are satisfied. The group $\Z^2$ can be bi-ordered,
say lexicographically. (In fact, there are uncountably many different
bi-orders on $\Z^2$; e.g.\ there are already two for each line of irrational
slope in $\R^2$ through $(0,0)$.)

All that remains to be proven is that $F \cong \pi_1(\widetilde{M})$ has a
bi-ordering which is invariant under conjugation by elements in
$G = \pi_1(M)$. We note that $F$ is an infinitely-generated free
group.  There is a free basis for $F$ consisting of
the generators $x_{i,j}$ represented by a loop that goes around the
central
curve of the M\"obius band $B_{i,j}$, connected by a tail to the
basepoint in some (non-canonical) way; for definiteness we may take
$$
x_{i,j} = a^ib^jcb^{-j}a^{-i} \ \ \ \ \ ((i, j)\in \Z^2)
$$
as a free generating set for $F$.
Now $G$ acts upon $F$ by conjugation, which may be described in terms
of the generators as follows.

\begin{lemma}\label{conj}\sl
Suppose $g \in G$ has exponent sums $m$ and $n$ in $a$ and $b$,
respectively.  Then $$gx_{i,j}g^{-1} = w_{i,j}x_{i+m,j+n}w_{i,j}^{-1},$$
where $w_{i,j}=w_{i,j}(g) \in F$.
\end{lemma}

\begin{proof}[Proof of the lemma]
We just take $w_{i,j} = ga^ib^{-n}a^{-i-m}$; by calculating the exponent
sums of $a$ and $b$ in $w_{i,j}$ we can verify that indeed $w_{i,j} \in F$.
\end{proof}

For the following, $F_{ab}$ denotes the abelianization of $F$, which
is an infinitely generated free abelian group, with generators, say
$\widetilde{x}_{i,j}$; the abelianization map $F \to F_{ab}$ is just
$x_{i,j} \to \widetilde{x}_{i,j}$.  Any automorphism $\phi$ of $F$
induces a unique automorphism $\phi_{ab}$ of $F_{ab}$.  For example, in
the above lemma, under abelianization the conjugation map is just the shift
$\widetilde{x}_{i,j} \to \widetilde{x}_{i+m,j+n}$. Now lemma \ref{shiftinv}
completes the proof of proposition \ref{3P2}.\end{proof}

\begin{lemma}\label{shiftinv}\sl
There is a bi-ordering of the free group
$F = \langle x_{i,j}\rangle; (i,j) \in \Z^2$ which is invariant under
every automorphism $F \to F$ which induces, on $F_{ab}$, a uniform shift
automorphism $\widetilde{x}_{i,j} \to \widetilde{x}_{i+m,j+n}$.
\end{lemma}

\begin{proof} We use the Magnus expansion \cite{MKS}, sending $F$ into the ring of
formal power series in the infinitely many noncommuting variables
$X_{i,j}$.
Since there are infinitely many generators, some care must be taken in
defining the ring $\Z[[X_{i,j}]]$, which we take to be the ring
consisting of formal power series in the $X_{i,j}$, but we consider only
such series which involve just finitely many different variables.  The
Magnus map $\mu\co F \to \Z[[X_{i,j}]]$ is given by
$$ \mu(x_{i,j}) = 1 + X_{i,j};
\quad   \mu(x_{i,j}^{-1}) = 1 - X_{i,j} + X_{i,j}^2 - X_{i,j}^3 + \cdots
$$
Clearly the image of $F$ lies in the group of units of the form
$\{1 + O(1)\}$ inside  $\Z[[X_{i,j}]]$, and it is an embedding of groups,
by the same proof as in \cite{MKS}.   Elements of  $\Z[[X_{i,j}]]$ may be
written in standard form, arranged in ascending degree, and within a
degree terms are arranged lexicographically by their subscripts (which in
turn are ordered lexicographically).  Then two series are compared
according to the
coefficient of the first term at which they differ (here is where the
finiteness assumption is necessary).  The proof that this defines a
(multiplicative)
bi-invariant ordering of $\{1 + O(1)\}$, is routine - cf.\ \cite{KR}.
Via the injection $\mu$, we may regard $F$ as a subgroup, and hence it
is also biordered.

Finally, we argue that this ordering has the desired invariance property;
equivalently, that $\phi$ preserves the positive cone of $F$.
Consider an automorphism $\phi\co F \to F$ such that
$\phi_{ab}\co F_{ab} \to F_{ab}$ is a shift
$\phi_{ab}(\widetilde{x}_{i,j}) = \widetilde{x}_{i+m,j+n}.$
This means that $\phi(x_{i,j}) = x_{i+m,j+n} c_{i,j}$, where $c_{i,j}$ is
in the commutator subgroup $[F,F].$  Since $[F,F]$ maps into
$\{1 + O(2)\}$ under the Magnus embedding, the effect of
$\phi$ is reflected in $\Z[[X_{i,j}]]$ by the substitution $X_{i,j} \to
X_{i+m,j+n} + O(2).$  Therefore, if the Magnus expansion of $w \in F$ is
$\mu(w) = 1 + w_d(X_{i,j}) + O(d+1)$, where $w_d$ is the sum of all degree
$d$ terms, then
$\mu\phi(w) = 1 + w_d(X_{i+m,j+n}) +$ terms of higher degree.
Therefore, the lowest degree nonzero terms of the Magnus expansions of $w$
and $\phi(w)$ are identical, except that the subscripts are shifted.  Thus
the ``first'' nonconstant terms of both $w$ and $\phi(w)$ have the same
coefficient, and we conclude that $\phi$ preserves the positive cone of
$F$ in the ordering we described. \end{proof}


\section{Punctured-torus bundles over the circle}

Suppose that a 3-manifold $M^3$ is a fibre bundle over the
circle, with fibre some surface $S$, and monodromy map $\f\co S \to S$.
A necessary and sufficient condition for $\pi_1(M^3)$ to be bi-ordered is
that there exists some bi-order on $\pi_1(S)$ which is invariant under
$\f$.
This can be seen by considering the short exact sequence of the fibration
$$
1 \longrightarrow \pi_1(S) \ \hookrightarrow \ \pi_1(M^3) \longrightarrow
\Z
\longrightarrow 1
$$
and applying lemma \ref{extensions}.
It is an important and apparently hard problem in general to decide which
automorphisms of free groups or surface groups leave some bi-ordering
invariant.

In the present paper we restrict our attention to the case where $S$ is
a once-punctured torus - for a different line of attack see \cite{PR}.
We shall apply a slightly more sophisticated version of
lemma \ref{shiftinv} in order to prove

\begin{thm}\sl \label{41bi}
Suppose the 3-manifold $M^3$ is a fibre bundle over the
circle, with fibre a once-punctured torus $S$. Suppose that the monodromy
map $\f\co S \to S$ is orientation preserving and induces in homology a
homomorphism $\Phi\co H_1(S) \to H_1(S)$ which preserves some bi-ordering
on
$H_1(S)\cong \Z^2$. Then $\pi_1(M^3)$ is bi-orderable.
\end{thm}

It is known which automorphisms of $\Z^n$ leave some bi-order invariant.
This is due to Levitt \cite{Levitt}: suppose $T\co \C^n \to \C^n$ is a
linear
transformation that restricts to an automorphism $T\co \Z^n \to \Z^n$
(i.e.\ it is represented by some $n\times n$ matrix with integer entries,
and determinant $\pm 1$). Consider a basis $\mathcal{B}$ of $\C^n$ such
that the matrix of $T$ with respect to $\mathcal{B}$ is in Jordan normal
form. We can split $\C^n$ uniquely as a direct sum $E\oplus F$, where we
define $E$ to be spanned by the
vectors of $\mathcal{B}$ belonging to Jordan blocks with positive real
eigenvalue, and $F$ to be spanned by the vectors of $\mathcal{B}$ with
negative or complex eigenvalues. Now Levitt's criterion is: the
homomorphims
$T\co \Z^n \to \Z^n$ leaves some bi-order on $\Z^n$ invariant if and only
if $F$ does not intersect the integer lattice $\Z^n \subseteq \C^n$ except
in $0$. The proof
is not difficult, and left to the reader - the basic observation is that
for any bi-ordering of $\Z^n$ there exists a hyperplane in $\R^n$ such
that all integer lattice points to one side of the plane are in the
positive, and those on the other side in the negative cone.

We remark that if the monodromy map is orientation reversing, then
$\pi_1(M^3)$ is definitely not bi-orderable, because it contains a Klein
bottle group. Also, if the monodromy map is periodic, then $\pi_1(M^3)$
cannot be bi-orderable, because a periodic automorphism of $\pi_1(S)$
cannot leave a bi-ordering of $\pi_1(S)$ invariant. This proves, for
example, that the group of the trefoil knot, which is a fibred knot with
a monodromy map of period $6$, is not bi-orderable. This fact, which was
pointed out in \cite{N2}, seems to have dissuaded mathematicians from
studying bi-orderability of knot groups for the last 25 years.

\medskip
{\bf Examples for theorem \ref{41bi}}\qua An orientation preserving
(i.e.\ determinant 1) automorphism $\Phi$ of $\Z^2$
preserves a bi-order if and only if its eigenvalues are
both positive (but not necessarily distinct). For if $\Phi$ has
eigenvalues
$\lambda_1, 1/\lambda_1>0$, where $\stackrel{\to}{v_1}$ is the eigenvector
for $\lambda_1$ and $\stackrel{\to}{v_2}$ is the other basis vector in a
Jordan normal basis, then we can define the positive cone of a bi-ordering
of $\Z^2$ by
$$
P_{\Z^2} = \{
\stackrel{\to}{z}=c_1\stackrel{\to}{v_1}+c_2\stackrel{\to}{v_2}
\ \in \Z^2 \ : \ c_2>0 {\rm \ or \ } (c_2=0 {\rm \ and \ } c_1>0) \ \}.
$$
We observe that this order is invariant under $\Phi$.

For instance,
theorem \ref{41bi} implies that the complement of the figure-of-eight knot
has bi-orderable fundamental group, because it is a fibred knot with
punctured torus fibre $S$, and the matrix of the monodromy
action on $H_1(S)$ is {\tiny${\pmatrix{2 & 1 \cr 1 & 1}}$}, which has 
two positive eigenvalues. The knot $4_1$ is the only \emph{classical}
knot (in $S^3$) which is covered by theorem \ref{41bi}, because
the only classical fibred knots of genus $1$ are the knots $3_1$
and $4_1$. A few more classical knots are bi-orderable by the main
theorem of \cite{PR}.

If, on the other hand, $\Phi$ has two negative, or two complex conjugate
eigenvalues, then it cannot respect any decomposition of $\Z^2$ into a
``positive'' and a ``negative'' half-space; thus $\Phi$ cannot preserve
any ordering on $\Z^2$.

{\bf Proof of theorem \ref{41bi}}\qua Consider the short exact sequence
$$
1 \longrightarrow \pi_1(\widetilde{S})\longrightarrow \pi_1(S)
\stackrel{ab}{\longrightarrow} H_1(S) \longrightarrow 1,
$$
where $ab$ is the abelianisation homomorphism (which geometrically
corresponds to ``patching the puncture in $S$''), and $\widetilde{S}$ is
the covering space of $S$ whose fundamental group is the commutator
subgroup
$[\pi_1(S), \pi_1(S)]$.
One can picture $\widetilde{S}$ as the plane $\R^2$ with all integer
lattice points $\Z^2$ removed.

In particular, $\pi_1(\widetilde{S})$ is an infinitely generated free
group,
with generating set
$\{x_{i,j} \ | \ (i,j)\in \Z^2 \}$, where the loop $x_{i,j}$ has winding
number
one around the puncture at the point $(i,j) \in \R^2$, and winding number
zero around the punctures at all the other integer lattice points (many
different choices are possible here). Note that the abelianisation
$H_1(\widetilde{S})$ is an infinitely generated free abelian group with
generators $\{\widetilde{x}_{i,j} : (i,j)\in \Z^2\}$.

We have to prove that $\phi_*\co \pi_1(S)\to \pi_1(S)$ leaves some
bi-order on $\pi_1(S)$ invariant. We already have a $\Phi$-invariant
bi-order on $H_1(S)$, and in fact we shall take it to be the order
constructed in the example for theorem \ref{41bi}. Thus
it suffices by lemma \ref{extensions} to find a
bi-order on $\pi_1(\widetilde{S})$ which is
\begin{itemize}
\item[(i)]invariant under conjugation by elements in $\pi_1(S)$, and
\item[(ii)]invariant under (the restriction of) $\phi_*$.
\end{itemize}

Let's study the effect of conjugation by elements in $\pi_1(S)$ and
of $\phi_*$. Conjugation by an element $g \in \pi_1(S)$ with $ab(g)=(n,m)
\in \Z^2$ sends a generator $x_{i,j}$ to $w_{i,j}x_{i+m,j+n}w_{i,j}^{-1}$,
where $w_{i,j} \in \pi_1(\widetilde{S})$
depends on $g$, as in lemma \ref{conj}. Thus, on the abelianisation
$H_1(\widetilde{S})$, conjugation by $g$ induces a uniform shift
automorphism $\widetilde{x}_{i,j} \to \widetilde{x}_{i+m,j+n}$. (Note that
the mapping $\Z^2 \to \Z^2$, $(i,j) \mapsto (i+m, j+n)$ preserves
our ordering of $\Z^2$.)

Moreover, the restriction of $\phi_*$ to $\pi_1(\widetilde{S})$
is given geometrically by the action of the lift $\widetilde{\phi}\co
\widetilde{S} \to \widetilde{S}$ which fixes the point $(0,0)$. This
sends the generator $x_{i,j}$ to a conjugate of the generator
$x_{\Phi(i,j)}$, and induces on the abelianisation $H_1(\widetilde{S})$
an automorphism determinded by another simple permutation of the
generators:
$\widetilde{x}_{i,j} \mapsto \widetilde{x}_{\Phi(i,j)}$. The proof of
theorem \ref{41bi} is now completed by the following

\begin{lemma}\label{geninv}\sl
There is a bi-ordering of the free group
$\pi_1(\widetilde{S}) = \langle x_{i,j}\rangle; (i,j) \in \Z^2$ which is
invariant under every automorphism which induces, on the abelianisation
$H_1(\widetilde{S})$, an automorphism that acts simply by permuting
the variables: $\widetilde{x}_{i,j} \to \widetilde{x}_{\sigma(i,j)}$,
where the permutation $\sigma\co \Z^2 \to \Z^2$ preserves the bi-ordering
of $\Z^2$.
\end{lemma}
\begin{proof} The proof is virtually the same as for lemma \ref{shiftinv}.
We order the variables $\{X_{i,j}\}$ by defining $X_{i,j}$ to ``come
before''
$X_{i',j'}$ if $(i,j)>(i',j')$ in the $\Phi$-invariant bi-order of $\Z^2$.
Then in order to define which of two words in $\pi_1(\widetilde{S})$ is
the larger we compare their images under the Magnus map $\mu$. Finally
the invariance property is proved precisely as in lemma
\ref{shiftinv}.\end{proof}

\Addresses\recd
\end{document}